# $\phi$-PARABOLICITY AND THE UNIQUENESS OF SPACELIKE HYPERSURFACES IMMERSED IN A SPATIALLY WEIGHTED GRW SPACETIME

ALMA L. ALBUJER[1,*], HENRIQUE F. DE LIMA[2], ARLANDSON M. OLIVEIRA[2]
AND MARCO ANTONIO L. VELÁSQUEZ[2]

ABSTRACT. In this paper, we extend a technique due to Romero, Rubio and Salamanca [24, 25, 26] establishing sufficient conditions to guarantee the parabolicity of complete spacelike hypersurfaces immersed in a weighted generalized Robertson-Walker spacetime whose fiber has $\phi$-parabolic universal Riemannian covering. As some applications of this criteria, we obtain uniqueness results concerning spacelikes hypersurfaces immersed in spatially weighted generalized Robertson-Walker spacetimes. Furthermore, Calabi-Bernstein type results are also given.

## 1. INTRODUCTION

The study of immersed spacelike hypersurfaces with constant mean curvature in a Lorentzian manifold has attracted the interest of a considerable group of geometers, as evidenced by the amount of works that it has generated in the last decades. This is due not only to its mathematical interest, but also to its relevance in General Relativity. For example, constant mean curvature spacelike hypersurfaces are particularly suitable for studying the propagation of gravitational radiation (cf. [28]; see also [20] for a summary of several reasons justifying this interest).

From the mathematical point of view, the study of the geometry of constant mean curvature spacelike hypersurfaces is mostly due to the fact that they exhibit nice Bernstein type properties. In this context, several authors have recently treated the problem of uniqueness for complete constant mean curvature spacelike hypersurfaces of generalized Robertson-Walker (GRW) spacetimes, that is, Lorentzian warped products with 1-dimensional negative definite base and Riemannian fiber. Along this branch, we point out the works of Romero, Rubio and Salamanca [24, 25, 26], where the authors studied non-compact complete spacelike hypersurfaces in GRW spacetimes whose fiber has a parabolic universal Riemannian covering. In this setting, they guaranteed the parabolicity of complete spacelike hypersurfaces, under boundedness assumptions on the warping function restricted to a spacelike hypersurface and on the hyperbolic angle function of the hypersurface. As an application of this new technique, they obtained several uniqueness results on complete maximal spacelike hypersurfaces.

Going a step further, our aim in this paper is to extend the previous mentioned parabolicity criteria to the context of complete constant weighted mean curvature spacelike hypersurfaces immersed in a weighted GRW spacetime. As some applications of this criteria, we obtain uniqueness results concerning these spacelike hypersurfaces in the case where the GRW spacetime is spatially weighted.

We recall that, given a pseudo-Riemannian manifold $(M, g)$ and a smooth function $\phi \in C^\infty(M)$, the weighted manifold $M_\phi$ associated to $M$ and $\phi$ is the triple $(M, g, d\mu = e^{-\phi}dM)$, where $dM$ is the standard volume element of $M$. In this context, Bakry and Émeri introduced in [7] the Bakry-Émeri-Ricci tensor as a suitable generalization of the classical Ricci tensor, Ric, defined by

$$\mathrm{Ric}_\phi = \mathrm{Ric} + \mathrm{Hess}\,\phi,$$

where Hess stands for the Hessian with respect to the metric $g$. So, it seems natural to try to extend results stated in terms of the Ricci curvature tensor to analogous results for the Bakry-Émeri-Ricci tensor. It is also interesting to remark that weighted manifolds are closely related to the study of Ricci solitons, harmonic heat flows, and many others concepts that have a great mathematical interest. For an overview of results in this context, we refer the reader to [27] and [29].

This manuscript is organized as follows: in Section 2, we introduce some basic notions and facts related to spacelike hypersurfaces in GRW spacetimes and weighted manifolds in general. In Section 3, we establish our $\phi$-parabolicity criteria concerning complete spacelike hypersurfaces immersed in a weighted GRW spacetime. Specifically, in Theorem 1 we prove that

> Let $S$ be a complete spacelike hypersurface in a weighted GRW spacetime $M_\phi$ with weight function $\phi$, whose fiber $F$ is complete with $\phi$-parabolic universal Riemannian covering. If the hyperbolic angle function $\Theta$ of $S$ is bounded and the restriction $f(h)$ on $S$ of the warping function $f$ of $M$ satisfies $\sup_S f(h) < \infty$ and $\inf_S f(h) > 0$, then $S$ is $\phi$-parabolic.

In Section 4, we apply the above criteria in order to get uniqueness results for spacelike hypersurfaces in spatially weighted GRW spacetimes. We also establish some nice consequences for the particular cases where $M$ is static (i.e. with constant warping function) or a steady state type spacetime (i.e. with warping function $f(t) = e^t$). Finally, in Section 5 we present new Calabi-Bernstein type results establishing non-parametric versions of the results obtained in Section 4.

## 2. Set up

Let $(F, g_F)$ be a connected, $n$-dimensional, oriented Riemannian manifold, $I \subset \mathbb{R}$ an open interval endowed with the metric $-dt^2$ and $f : I \to \mathbb{R}$ a positive smooth function. Let us consider the product manifold $M = I \times F$. The class of Lorentzian manifolds which will be of our concern here is the one obtained by furnishing $M$ with the Lorentzian metric

$$(2.1) \qquad \overline{g} = -\pi_I^*(dt^2) + f^2(\pi_I)\pi_F^*(g_F),$$

where $\pi_I$ and $\pi_F$ are the projections onto the factors $I$ and $F$, respectively. In this case, we simply write

$$(2.2) \qquad M = -I \times_f F.$$

Let us observe that $M$ is a Lorentzian warped product with warping function $f$ and fiber $F$.

When $F$ has constant sectional curvature, the warped product (2.2) has been known in the mathematical literature as a Robertson-Walker (RW) spacetime, an allusion to the fact that, for $n = 3$, it is an exact solution of the Einstein's field equations (cf. Chapter 12 of [23]). In the general case, after [6], the warped product (2.2) has usually been referred to as a *generalized Robertson-Walker (GRW) spacetime*, and we will stick to this usage along this paper.

Let $S$ be an $n$-dimensional connected manifold. A smooth immersion $x : S \to M$ is said to be a *spacelike hypersurface* if $S$, furnished with the metric induced from (2.1) via $x$, is a Riemannian manifold. If this is so, we will always assume that the metric on $S$ is the induced one, which will be denoted by $g_S$. In this setting, it follows from the connectedness of $S$ that one can uniquely choose a globally defined timelike normal unit vector field $N \in \mathfrak{X}(S)^\perp$, having the same time-orientation of $\partial_t$, i.e., such that $\overline{g}(N, \partial_t) \leq -1$, where $\partial_t$ is the coordinate vector field induced by the universal time on $M$. One then says that $N$ is the *future-pointing Gauss map* of $S$.

We will consider two particular functions naturally attached to a spacelike hypersurface $S$ immersed into a GRW spacetime $M$, namely, the *(vertical) height function* $h = (\pi_I)|_S$ and the *hyperbolic angle function* $\Theta = \overline{g}(N, \partial_t)$. Given any vector field $V \in \mathfrak{X}(M)$, we denote by $V^F = (\pi_F)_*(V) = V + \overline{g}(V, \partial_t)\partial_t$ the projection of $V$ onto $F$. In particular, $N^F = N + \Theta\partial_t$, and therefore

$$(2.3) \qquad \left|N^F\right|^2 = \overline{g}(N^F, N^F) = \Theta^2 - 1.$$

Let us denote by $\overline{\nabla}$ and $\nabla$ the gradients with respect to the metrics of $M$ and $S$, respectively. Then, a simple computation shows that the gradient of $\pi_I$ on $M$ is given by

$$\overline{\nabla}\pi_I = -\overline{g}(\overline{\nabla}\pi_I, \partial_t)\partial_t = -\partial_t,$$

so that the gradient of $h$ on $S$ is

$$(2.4) \qquad \nabla h = (\overline{\nabla}\pi_I)^\top = -\partial_t^\top = -\partial_t - \Theta N.$$

Thus, from (2.4) we get

$$(2.5) \qquad |\nabla h|^2 = \Theta^2 - 1,$$

where $|\nabla h|^2 = g_S(\nabla h, \nabla h)$ denotes the norm of the vector field $\nabla h$ on $S$.

At this point we recall that, given a semi-Riemannian manifold $(P, g)$ and a smooth function $\phi$ on $P$, the *weighted manifold* $P_\phi$ associated to $P$ and $\phi$ is the triple $(P, g, d\mu = e^{-\phi}dP)$, where $dP$ is the canonical volume element of $P$. In this setting, we will consider the so-called *Bakry-Émery-Ricci tensor*, introduced by Bakry and Émery in [7] as a suitable extension of the standard Ricci tensor Ric, which is defined by

$$\text{Ric}_\phi = \text{Ric} + \text{Hess}\,\phi. \tag{2.6}$$

For a spacelike hypersurface $S$ immersed in a weighted GRW spacetime $M$, the $\phi$-*divergence operator* on $S$ is defined by

$$\text{div}_\phi(X) = e^\phi \text{div}(e^{-\phi}X),$$

for any tangent vector field $X$ on $S$ and, given a smooth function $u : S \to \mathbb{R}$, its *drifted Laplacian* is defined by

$$\Delta_\phi u = \text{div}_\phi(\nabla u) = \Delta u - g_S(\nabla u, \nabla \phi). \tag{2.7}$$

Furthermore, according to Gromov [14], the $\phi$-*mean curvature* $H_\phi$ of $S$ is defined by

$$nH_\phi = nH - \overline{g}(\overline{\nabla}\phi, N), \tag{2.8}$$

where $H$ denotes the standard mean curvature of $S$ with respect to its future-pointing Gauss map $N$.

## 3. $\phi$-Parabolicity of spacelike hypersurfaces

A smooth function $u$ on a weighted manifold $P_\phi$ is said to be $\phi$-*superharmonic* if $\Delta_\phi u \leq 0$. Taking this into account, the weighted manifold $(P, g, d\mu = e^{-\phi}dP)$ is called $\phi$-*parabolic* if there is no nonconstant, nonnegative, $\phi$-superharmonic function on $P$. On the other hand, for any compact subset $K \subset P$, we define the $\phi$-*capacity* of $K$ as

$$\text{cap}_\phi(K) = \inf\left\{\int_P |\nabla u|^2 d\mu : u \in \text{Lip}_0(P) \text{ and } u\big|_K \equiv 1\right\},$$

where $\text{Lip}_0(P)$ is the set of all compactly supported Lipschitz functions on $P$. The following statement relates the notion of $\phi$-capacity to the concept of $\phi$-parabolicity (cf. [13, Proposition 2.1]).

**Lemma 1.** *The weighted manifold $(P, g, d\mu = e^{-\phi}dP)$ is $\phi$-parabolic if, and only if, $\text{cap}_\phi(K) = 0$ for any compact set $K \subset P$.*

Let us recall that given two Riemannian manifolds $(P, g)$ and $(P', g')$, a diffeomorphism $\psi$ from $P$ onto $P'$ is called a *quasi-isometry* if there exists a constant $c \geq 1$ such that

$$c^{-1}|v|_g \leq |d\psi(v)|_{g'} \leq c|v|_g$$

for all $v \in T_pP, p \in P$ (cf. [18] for more details). Suppose that we can endow both $P$ and $P'$ with the same weight function $\phi$. We can reason as in Section 5 of [12] to verify that the $\phi$-capacity changes under a quasi-isometry at most by a constant factor. So, from Lemma 1, we can state the following result.

**Lemma 2.** *Let $(P, g)$ and $(P', g')$ be two Riemannian manifolds endowed with the same weight function $\phi$. If $P$ and $P'$ are quasi-isometric, then $P$ and $P'$ are $\phi$-parabolic or not simultaneously.*

Given a spacelike hypersurface $x : S \to M$ in a GRW spacetime $M = -I \times_f F$, the following lemma provides sufficient conditions to guarantee that the hypersurface $S$ and the fiber $F$ are quasi-isometric (cf. [24, Lemma 4.1]).

**Lemma 3.** *Let $x : S \to M$ be a spacelike hypersurface in a GRW spacetime $M = -I \times_f F$, whose hyperbolic angle function $\Theta$ is bounded. If the warping function $f$ on $S$ satisfies*
   (i) $\sup_S f(h) < \infty$ *and*
   (ii) $\inf_S f(h) > 0$,
*then $\pi = \pi_F \circ x$ is a quasi-isometry from $S$ onto $F$.*

We can now present the main result of this section.

**Theorem 1.** *Let $S$ be a complete spacelike hypersurface in a weighted GRW spacetime $M_\phi$ with weight function $\phi$, whose fiber $F$ is complete with $\phi$-parabolic universal Riemannian covering. If the hyperbolic angle function $\Theta$ of $S$ is bounded and the restriction $f(h)$ on $S$ of the warping function $f$ of $M$ satisfies*

(i) $\sup_S f(h) < \infty$ *and*
(ii) $\inf_S f(h) > 0$,

*then $S$ is $\phi$-parabolic.*

The proof of Theorem 1 follows the same steps of [24, Theorem 4.4] (see also [25, Theorem 1]). For this, we will need a standard result on covering spaces (cf. [15] for instance).

**Lemma 4.** *Let $\rho : (\tilde{E}, \tilde{x}_0) \to (E, x_0)$ be a covering space and let $h : (W, y_0) \to (E, x_0)$ be a continuous map, where $W$ is a path connected and locally path connected topological space. Then, there exists a lift $\tilde{h} : (W, y_0) \to (\tilde{E}, \tilde{x}_0)$ of $h$ if, and only if, $h_*(\pi_1(W, y_0)) \subset \rho_*(\pi_1(\tilde{E}, \tilde{x}_0))$.*

*Proof of Theorem 1.* Using [6, Lemma 3.1], we know that the projection on the fiber, $\pi : S \to F$, is a covering map. Moreover, by Lemma 3, we can find a constant $c \geq 1$ such that

$$(3.1) \qquad c^{-1} g_F(d\pi(v), d\pi(v)) \leq g_S(v,v) \leq c g_F(d\pi(v), d\pi(v)),$$

for all $v \in T_p S$ and $p \in S$.

Let $(\tilde{S}, g_{\tilde{S}})$ be the universal Riemannian covering of $(S, g_S)$ and denote by $\tilde{\pi}_S : \tilde{S} \to S$ the corresponding Riemannian covering map. From Lemma 4 we conclude that there exists a lift $\tilde{h} : \tilde{S} \to \tilde{F}$ of the map $h = \pi \circ \tilde{\pi}_S : \tilde{S} \to F$. It is easy to check that $\tilde{h}$ is, in fact, a diffeomorphism from $\tilde{S}$ to $\tilde{F}$. Note that (3.1) implies

$$c^{-1} g_{\tilde{F}}(d\tilde{h}(\tilde{v}), d\tilde{h}(\tilde{v})) \leq g_{\tilde{S}}(\tilde{v}, \tilde{v}) \leq c g_{\tilde{F}}(d\tilde{h}(\tilde{v}), d\tilde{h}(\tilde{v})),$$

for any $\tilde{v} \in T_{\tilde{p}}\tilde{S}, \tilde{p} \in \tilde{S}$, which means that $\tilde{h}$ is a quasi-isometry from $(\tilde{S}, g_{\tilde{S}})$ onto $(\tilde{F}, g_{\tilde{F}})$.

Finally, let $u$ be a nonnegative $\phi$-superharmonic function on $S$, and put $\tilde{u} = u \circ \tilde{\pi}_S$. The function $\tilde{u}$ is a nonnegative $\phi$-superharmonic function on the $\phi$-parabolic Riemannian manifold $\tilde{S}$. Therefore, $\tilde{u}$ must be constant and, consequently, $u$ is also constant. □

As a direct consequence of Theorem 1 we get the following corollaries.

**Corollary 1.** *Let $S$ be a complete spacelike hypersurface in a weighted GRW spacetime $M_\phi$ with weight function $\phi$ and whose fiber $F$ is complete, simply connected and $\phi$-parabolic. If the hyperbolic angle function $\Theta$ of $S$ is bounded and the warping function on $S$, $f(h)$, is bounded and it satisfies $\inf_S f(h) > 0$, then $S$ is $\phi$-parabolic.*

We recall that a GRW is said to be *static* when its warping function is constant, which, without loss of generality, can be supposed equal to 1.

**Corollary 2.** *Let $S$ be a complete spacelike hypersurface in a static weighted GRW spacetime $M_\phi$ with weight function $\phi$ and whose fiber $F$ is complete with $\phi$-parabolic universal Riemannian covering. If the hyperbolic angle function $\Theta$ of $S$ is bounded, then $S$ is $\phi$-parabolic.*

## 4. Uniqueness results

It follows from a splitting theorem due to Case (cf. [8], Theorem 1.2) that if $M$ is a weighted GRW spacetime endowed with a weight function $\phi$ which is bounded and such that $\overline{\mathrm{Ric}}_\phi(V, V) \geq 0$ for any timelike vector field $V$, then $\phi$ must be constant along $\mathbb{R}$. Motivated by this result, in what follows we will deal with *spatially weighted GRW spacetimes* $M$, which means that the weight function $\phi$ does not depend on the parameter $t \in I$, that is, $\overline{g}(\overline{\nabla}\phi, \partial_t) = 0$.

**Remark 1.** *We note that the $\phi$-mean curvature of a slice $\{t_0\} \times F$ of a spatially weighted GRW spacetime is given by*

$$H_\phi(t_0) = (\log f)'(t_0).$$

*Indeed, since $\partial_t$ is a timelike normal vector field to the slices $\{t\} \times F$, from (2.8) we have that $H_\phi(t_0) = H(t_0) = (\log f)'(t_0)$.*

As an application of Theorem 1, we will prove in this section some uniqueness results concerning spacelike hypersurfaces immersed in a spatially weighted GRW spacetime. For this, we will also need the following lemma.

**Lemma 5.** *Let $S$ be a spacelike hypersurface immersed in a spatially weighted GRW spacetime $M_\phi$ with weight function $\phi$. Then,*

$$\Delta_\phi h = -(\ln f)'(h)(n + |\nabla h|^2) - n\Theta H_\phi, \tag{4.1}$$

$$\Delta_\phi \mathcal{F}(h) = -n(f'(h) + f(h)\Theta H_\phi), \tag{4.2}$$

*where $\mathcal{F}(t) = \int_{t_0}^t f(s)ds$, and*

$$\begin{aligned}\Delta_\phi(f(h)\Theta) =& nf(h)\overline{g}(\nabla H_\phi, \partial_t) + nf'(h)H_\phi + f(h)\Theta|A|^2 \\ & + f(h)\Theta\overline{\mathrm{Hess}}\,\phi(N,N) \\ & + f(h)\Theta\left(\mathrm{Ric}^F(N^F, N^F) - (n-1)(\log f)''(h)|\nabla h|^2\right),\end{aligned} \tag{4.3}$$

*where $\mathrm{Ric}^F$ stands for the Ricci curvature tensor of the fiber $F$.*

*Proof.* Equations (4.1) and (4.2) correspond to [9, Lemma 1], and (4.3) corresponds to [4, Lemma 1]. However, until the moment reference [4] is not yet published, so for the sake of completeness we present it here a proof of (4.3).

In [5, Corollary 8.2] it is proven that

$$\begin{aligned}\Delta(f(h)\Theta) =& nf(h)\overline{g}(\nabla H, \partial_t) + nf'(h)H + f(h)\Theta|A|^2 \\ & + f(h)\Theta\left(\mathrm{Ric}^F(N^F, N^F) - (n-1)(\log f)''(h)|\nabla h|^2\right).\end{aligned} \tag{4.4}$$

Taking into account (2.8), it follows that

$$nf(h)\overline{g}(\nabla H, \partial_t) = nf(h)\overline{g}(\nabla H_\phi, \partial_t) + f(h)\partial_t^\top \overline{g}(\overline{\nabla}\phi, N).$$

Moreover, from a straightforward computation we get

$$\partial_t^\top\left(\overline{g}\left(\overline{\nabla}\phi, N\right)\right) = -\frac{f'}{f}(h)\overline{g}(\overline{\nabla}\phi, N) + \Theta\,\overline{\mathrm{Hess}}\,\phi(N,N) - \overline{g}(\overline{\nabla}\phi, A\partial_t^\top),$$

and $\nabla(f(h)\Theta) = -f(h)A\partial_t^\top$. So (4.4) can be written as

$$\begin{aligned}\Delta(f(h)\Theta) =& nf(h)\overline{g}(\nabla H_\phi, \partial_t) - f'(h)\overline{g}(\overline{\nabla}\phi, N) + f(h)\Theta\overline{\mathrm{Hess}}\,\phi(N,N) \\ & + \overline{g}(\overline{\nabla}\phi, \nabla(f(h)\Theta)) + nf'(h)H + f(h)\Theta|A|^2 \\ & + f(h)\Theta\left(\mathrm{Ric}^F(N^F, N^F) - (n-1)(\log f)''(h)|\nabla h|^2\right).\end{aligned} \tag{4.5}$$

Finally, (4.3) follows from (4.5) and (2.7). $\square$

We can now prove our first uniqueness result. A spacelike hypersurface is called $\phi$-*maximal* if its $\phi$-mean curvature $H_\phi$ is identically zero. A slab $[t_1, t_2] \times F = \{(t,q) \in M : t_1 \leq t \leq t_2\}$ of a GRW spacetime $M$ is called a *timelike bounded region*.

**Theorem 2.** *Let $M_\phi$ be a spatially weighted GRW spacetime whose fiber $F$ is complete with $\phi$-parabolic universal Riemannian covering, and such that the warping function $f$ is monotone. The only $\phi$-maximal complete spacelike hypersurfaces contained in a timelike bounded region of $M$ and with bounded hyperbolic angle function $\Theta$ are the slices $\{t_0\} \times F$, where $t_0 \in I$ is such that $f'(t_0) = 0$.*

*Proof.* Let $S$ be such a spacelike hypersurface. From Lemma 5, we obtain

$$\Delta_\phi \mathcal{F}(h) = -nf'(h).$$

Consequently, the monotonicity of $f$ implies that $\Delta_\phi \mathcal{F}(h)$ is globally either nonpositive or nonnegative signed. Since $S$ is contained in a timelike bounded region of $M$ and the warping function $f$ is monotone, the function $\mathcal{F}(h)$ is clearly bounded on $S$. From Theorem 1 we know that $S$ is $\phi$-parabolic, so $\mathcal{F}(h)$ is constant in $S$ and, hence, $h$ must be also constant in $S$. $\square$

In our next result, we will assume that the ambient space obeys the so-called *null convergence condition* (NCC). We recall that a GRW spacetime $M$ satisfies NCC if

$$\text{Ric}^F \geq (n-1)f^2(\log f)'' g_F, \tag{4.6}$$

which is equivalent to the Ricci curvature of $M$ being nonnegative on null or lightlike directions (cf. [21]).

**Theorem 3.** *Let $M_\phi$ be a spatially weighted GRW spacetime satisfying (4.6), with convex weight function $\phi$ (that is, $\overline{\text{Hess}}\,\phi \geq 0$) and whose fiber $F$ is complete with $\phi$-parabolic universal Riemannian covering. Let $S$ be a complete $\phi$-maximal spacelike hypersurface immersed in $M$, with bounded hyperbolic angle function $\Theta$ and such that the restriction $f(h)$ on $S$ of the warping function $f$ of $M$ satisfies*
  (i) $\sup f(h) < \infty$ *and*
  (ii) $\inf f(h) > 0$.
*Then $S$ is totally geodesic. In addition, if the inequality (4.6) is strict for all non-zero vector fields on $F$ or $\phi$ is strictly convex on $F$, then $S$ is a slice $\{t_0\} \times F$, where $t_0 \in I$ is such that $f'(t_0) = 0$.*

*Proof.* From Lemma 5, we have that the drifted Laplacian of the bounded function $f\Theta$ is given by

$$\begin{aligned}\Delta_\phi(f(h)\Theta) =\,& f(h)\Theta|A|^2 + f(h)\Theta\overline{\text{Hess}}\,\phi(N,N) \\ & + f(h)\Theta\left(\text{Ric}^F(N^F, N^F) - (n-1)(\log f)''(h)|\nabla h|^2\right).\end{aligned} \tag{4.7}$$

Since $\phi$ is convex and we are assuming that the null convergence condition (4.6) holds, it follows that $\Delta_\phi(f(h)\Theta) \leq 0$. Theorem 1 assures that $S$ is $\phi$-parabolic, so $f(h)\Theta$ must be constant. Therefore, returning to (4.7), we infer that $|A| \equiv 0$, that is, $S$ is totally geodesic,

$$\overline{\text{Hess}}\phi(N,N) = \text{Hess}^F\phi(N^F, N^F) = 0, \tag{4.8}$$

and

$$\text{Ric}^F(N^F, N^F) - (n-1)(\log f)''(h)|\nabla h|^2 = 0.$$

Consequently, if the inequality (4.6) is strict, or if $\phi$ is strictly convex on $F$, then (4.7) also gives that $|N^F| = |\nabla h| = 0$ on $S$, that is, $S$ is a slice. □

When the ambient space is static, we obtain the following,

**Theorem 4.** *Let $M_\phi$ be a static spatially weighted GRW spacetime, whose fiber $F$ is complete with $\phi$-parabolic universal Riemannian covering and such that its Bakry-Émery-Ricci tensor $\text{Ric}^F_\phi$ is nonnegative. Let $S$ be a complete spacelike hypersurface immersed in $M$ with constant $\phi$-mean curvature $H_\phi$. If the hyperbolic angle function $\Theta$ of $S$ is bounded, then $S$ is totally geodesic. In addition, if $\text{Ric}^F_\phi$ is definite positive at some point of $S$, then $S$ is a slice $\{t\} \times F$.*

*Proof.* From Lemma 5, (2.6), and (4.8) we get that

$$\Delta_\phi \Theta = (\text{Ric}^F_\phi(N^F, N^F) + |A|^2)\Theta. \tag{4.9}$$

Consequently, since we are supposing that $\text{Ric}^F_\phi$ is nonnegative and that $\Theta$ is negative and bounded on $S$, we can apply Corollary 2 to conclude that $\Theta$ is constant on $S$. Thus, returning to (4.9) we get that $|A| \equiv 0$, that is, $S$ is totally geodesic. Moreover, if $\text{Ric}^F_\phi$ is definite positive at some $p \in F$, considering once more equation (4.9) and taking into account relation (2.3), we conclude that $\Theta \equiv -1$ on $S$, which means that $S$ is a slice of $M$. □

We recall that the Gaussian space $\mathbb{G}^n$ corresponds to the Euclidean space $\mathbb{R}^n$ endowed with the Gaussian probability density $e^{-\phi(x)} = (2\pi)^{-\frac{n}{2}} e^{-\frac{|x|^2}{2}}$. From Corollary 3 of [16] we have that $\mathbb{G}^n$ has finite $\phi$-volume. Consequently, taking into account Remark 3.8 of [17], we conclude that $\mathbb{G}^n$ is $\phi$-parabolic.

On the other hand, we also recall that Xin [30] and Aiyama [1] proved simultaneously and independently that the only spacelike hypersurfaces of the Lorentz-Minkowski space $\mathbb{L}^{n+1} = -\mathbb{R} \times \mathbb{R}^n$, with constant mean curvature and having bounded hyperbolic angle function, are the spacelike hyperplanes. Hence, from Theorem 4 we get the following extension of this Xin-Aiyama result

**Corollary 3.** *The only complete spacelike hypersurfaces of $-\mathbb{R} \times \mathbb{G}^n$, with constant $\phi$-mean curvature and having bounded hyperbolic angle function, are the spacelike hyperplanes $\{t\} \times \mathbb{G}^n$.*

Proceeding, we will use a Bochner's formula due to Wei and Wylie [29] to obtain the following theorem.

**Theorem 5.** *Let $M_\phi$ be a static spatially weighted GRW spacetime endowed with a convex weight function $\phi$ and whose fiber $F$ is complete, with nonnegative sectional curvature, and such that its universal Riemannian covering is $\phi$-parabolic. Let $S$ be a complete spacelike hypersurface lying in a semi-space of $M$ and with constant $\phi$-mean curvature $H_\phi$. If the hyperbolic angle function $\Theta$ is bounded, then $S$ is a slice $\{t\} \times F$.*

*Proof.* Since we are supposing that $F$ has nonnegative sectional curvature, it follows from inequalities (3.3) and (3.4) of [11] that

$$\text{Ric}(X, X) \geq nH g_S(AX, X) + |AX|^2. \tag{4.10}$$

On the other hand, taking into account that

$$\text{Hess}\,\phi(X, X) = \overline{\text{Hess}}\,\phi(X, X) - \bar{g}(\overline{\nabla}\phi, N)g_S(AX, X),$$

from the convexity of the weight function $\phi$ we get

$$\text{Hess}\,\phi(X, X) \geq -\bar{g}(\overline{\nabla}\phi, N)g_S(AX, X). \tag{4.11}$$

From (4.10) and (4.11) we get the following lower bound for $\text{Ric}_\phi$

$$\text{Ric}_\phi(X, X) \geq nH_\phi g_S(AX, X) + |AX|^2. \tag{4.12}$$

Inequality (4.12) provides us

$$\text{Ric}_\phi(\nabla h, \nabla h) \geq nH_\phi g_S(A(\nabla h), \nabla h) + |A(\nabla h)|^2. \tag{4.13}$$

Since $H_\phi$ is constant, we have

$$\nabla \Delta_\phi h = -nH_\phi A(\nabla h). \tag{4.14}$$

On the other hand, from Bochner's formula (cf. [29])

$$\frac{1}{2}\Delta_\phi |\nabla h|^2 = |\text{Hess}\,h|^2 + g_S(\nabla h, \nabla \Delta_\phi h) + \text{Ric}_\phi(\nabla h, \nabla h). \tag{4.15}$$

Consequently, from (2.5), (4.13), (4.14) and (4.15) we get

$$\frac{1}{2}\Delta_\phi \Theta^2 = \frac{1}{2}\Delta_\phi |\nabla h|^2 \geq |\text{Hess}\,h|^2 \geq 0. \tag{4.16}$$

Thus, from Corollary 2 we have that $\Theta$ is constant and, returning to (4.16), we get $|\text{Hess}\,h|^2 = 0$ in $S$. Then, since $n|\text{Hess}\,h|^2 \geq (\Delta h)^2$, we have that $h$ is harmonic. So, since $\Delta h = -nH\Theta$, we also get that $H = 0$ in $S$ and from (4.10) we have that $S$ has nonnegative Ricci curvature. Therefore, since we are supposing that $S$ lies in a semi-space of $M$, we can apply the strong Liouville property due to Yau in [31] (see also [19, Theorem 4.8]) to conclude that $h$ must be constant and, hence, $S$ is a slice of $M$. □

Extending the ideas of [2], we now consider *spatially weighted steady state type spacetimes*, that is, GRW spacetimes of the type $-\mathbb{R} \times_{e^t} F$ whose fiber $F$ is endowed with a weight function $\phi$. In this setting, our next result is an extension of Theorem 8 of [2].

**Theorem 6.** *Let $M_\phi$ be a spatially weighted steady state type spacetime, whose fiber $F$ is complete, with nonnegative sectional curvature, and let $S$ be a complete spacelike hypersurface which lies in a timelike bounded region of $M$. Suppose that $|\overline{\nabla}\phi|$ is bounded on $S$. If the $\phi$-mean curvature $H_\phi$ of $S$ is constant and the hyperbolic angle function $\Theta$ is bounded, then $H_\phi = 1$. In addition, if the universal Riemannian covering of $F$ is $\phi$-parabolic, then $S$ is a slice $\{t\} \times F$.*

*Proof.* We claim that the mean curvature $H$ of $S$ is bounded. Indeed, since $\bar{g}(\overline{\nabla}\phi, \partial_t) = 0$, from (2.8) we have that

$$\begin{aligned} n|H| &\leq n|H_\phi| + |\bar{g}(\overline{\nabla}\phi, N)| \\ &= n|H_\phi| + |\bar{g}(\overline{\nabla}\phi, N^F)|. \end{aligned} \tag{4.17}$$

Thus, from (2.3) and (4.17) we get

$$n|H| \leq n|H_\phi| + |\overline{\nabla}\phi|(\Theta^2 - 1). \tag{4.18}$$

Consequently, since $H_\phi$ is constant and $\overline{\nabla}\phi$ and $\Theta$ are supposed to be bounded on $\Sigma^n$, it follows from (4.18) that $H$ is also bounded on $S$.

On the other hand, from inequality (16) of [2] we have that the Ricci curvature of $S$ satisfies

$$\mathrm{Ric}(X,X) \geq n - 1 - \frac{n^2 H^2}{4}.$$

So, we conclude that the Ricci curvature of $S$ is bounded from below.

Thus, we can apply the generalized maximum principle of Omori [22] and Yau [31] to guarantee that there exists a sequence $\{p_k\}$ in $S$ such that

$$\lim_k h(p_k) = \sup_S h, \quad \lim_k |\nabla h(p_k)| = 0 \text{ and } \limsup_k \Delta h(p_k) \leq 0.$$

But, from (2.5) and (2.7), we also get that

$$\lim_k \Theta(p_k) = -1 \text{ and } \limsup_k \Delta_\phi h(p_k) \leq 0.$$

Hence, taking into account formula (4.1), we can reason in a similar way to the proof of Theorem 8 of [2] to obtain that $H_\phi = 1$. Furthermore, from formula (4.2) we also have that

$$\Delta_\phi e^h = -ne^h(1+\Theta) \geq 0.$$

Therefore, assuming that the universal Riemannian covering of $F$ is $\phi$-parabolic, we can apply Theorem 1 to conclude that $h$ is constant in $S$. □

Our next result extends Theorem 5.3 of [10].

**Theorem 7.** *Let $M_\phi$ be a spatially weighted steady state type spacetime whose fiber $F$ is complete, with $\phi$-parabolic universal Riemannian covering, and let $S$ be a complete spacelike hypersurface immersed in $M$ with $H_\phi \geq 1$. If the hyperbolic angle function $\Theta$ of $S$ satisfies $-\Theta \leq H_\phi$, then $S$ is a slice $\{t\} \times F$.*

*Proof.* From formula (4.1) we get

$$\begin{aligned}
\Delta_\phi e^{-h} &= e^{-h}\left(|\nabla h|^2 - \Delta_\phi h\right) \\
&\leq ne^{-h}\left(|\nabla h|^2 + 1 + H_\phi \Theta\right) \\
&= ne^{-h}\Theta\left(\Theta + H_\phi\right).
\end{aligned} \tag{4.19}$$

Hence, our hypothesis on $\Theta$ guarantees that the function $e^{-h}$ is a $\phi$-superharmonic positive function on $S$. Therefore, we can apply once more Theorem 1 to get that $h$ is constant on $S$. □

## 5. Calabi-Bernstein type results

Let $\Omega \subseteq F$ be a connected domain and let $u \in \mathcal{C}^\infty(\Omega)$ be a smooth function, then $S(u)$ will denote the vertical graph over $\Omega$ determined by $u$, that is,

$$S(u) = \{(u(x), x) : x \in \Omega\} \subset M = -I \times_f F.$$

The graph is said to be entire if $\Omega = F$. The metric induced on $\Omega$ from the Lorentzian metric of the ambient space via $S(u)$ is

$$g_{S(u)} = -du^2 + f^2(u)g_F. \tag{5.1}$$

It can be easily seen that a graph $S(u)$ is a spacelike hypersurface if, and only if, $|Du|_F^2 < f^2(u)$, $Du$ being the gradient of $u$ in $F$ and $|Du|_F$ its norm, both with respect to the metric $g_F$. It is well known (cf. [6, Lemma 3.1]) that in the case where $F$ is a simply connected manifold, every complete spacelike hypersurface $S$ immersed in $M$ such that the warping function $f$ is bounded on $S$ is an entire spacelike graph over $F$. In particular, this happens for complete spacelike hypersurfaces contained in a timelike bounded region of $M$. It is interesting to observe that, in contrast to the case of graphs into a Riemannian space, an entire spacelike graph $S(u)$ in a Lorentzian spacetime is not necessarily complete, in the sense that the induced Riemannian metric is not necessarily complete on $F$. However, it can be proved that if $F$ is complete and $|Du|_F^2 \leq f^2(u) - c$ for certain positive constant $c > 0$, then $S(u)$ is complete. Although a particular case of this claim is proven in [3, Theorem 4.1], and the general proof is given in ([4, Proposition 9]), we will expose it here for the sake of completeness.

**Proposition 1.** *Let $F$ be a complete Riemannian manifold and $S(u)$ an entire spacelike vertical graph in $M = -I \times_f F$. If*
$$|Du|_F^2 \leq f^2(u) - c$$
*for certain positive constant $c > 0$, then $\Sigma^n(u)$ is complete.*

*Proof.* From (5.1), the Cauchy-Schwarz inequality and the assumptions of the proposition we get
$$g_{S(u)}(X,X) = -g_F(Du, X)^2 + f^2(u)g_F(X,X) \geq \left(f^2(u) - |Du|_F^2\right) g_F(X,X) \geq c g_F(X,X),$$
for every $X \in \mathfrak{X}(S(u))$. This implies that $L \geq \sqrt{c} L_F$, where $L$ and $L_F$ denote the length of a curve on $S(u)$ with respect to the Riemannian metrics $g_{S(u)}$ and $g_F$, respectively. As a consequence, as $F$ is complete by assumption, the induced metric on $S(u)$ from the metric of $M$ is also complete. □

The future-pointing Gauss map of a spacelike vertical graph $S(u)$ over $\Omega$ is given by the vector field

$$(5.2) \qquad N(x) = \frac{f(u(x))}{\sqrt{f^2(u(x)) - |Du(x)|_F^2}} \left( \partial_t|_{(u(x),x)} + \frac{1}{f^2(u(x))} Du(x) \right), \quad x \in \Omega.$$

Moreover, the shape operator $A$ of $S(u)$ with respect to its orientation (5.2) is given by

$$
\begin{aligned}
AX &= -\frac{1}{f(u)\sqrt{f^2(u) - |Du|_F^2}} D_X Du - \frac{f'(u)}{\sqrt{f^2(u) - |Du|_F^2}} X \\
&\quad + \left( \frac{-g_F(D_X Du, Du)}{f(u)(f^2(u) - |Du|_F^2)^{3/2}} + \frac{f'(u)g_F(Du, X)}{(f^2(u) - |Du|_F^2)^{3/2}} \right) Du,
\end{aligned}
\tag{5.3}
$$

for any tangent vector field $X$ tangent to $\Omega$. Consequently, if $S(u)$ is a spacelike vertical graph over a domain $\Omega$ of the fiber $F$ of a spatially weighted GRW spacetime $M$ endowed with a weight function $\phi$, it is not difficult to verify from (2.8) and (5.3) that the $\phi$-mean curvature function $H_\phi(u)$ of $S(u)$ is given by

$$nH_\phi(u) = -\mathrm{div}_\phi \left( \frac{Du}{f(u)\sqrt{f(u)^2 - |Du|_F^2}} \right) - \frac{f'(u)}{\sqrt{f(u)^2 - |Du|_F^2}} \left( n + \frac{|Du|_F^2}{f(u)^2} \right).$$

The differential equation $H_\phi(u) = 0$ with the constraint $|Du|_F < f(u)$ is called the *$\phi$-maximal spacelike hypersurface equation* in $M$, and its solutions provide $\phi$-maximal spacelike graphs in $M$.

From Theorem 2 we obtain the following Calabi-Bernstein type result.

**Theorem 8.** *Let $M_\phi$ be a spatially weighted GRW spacetime, whose fiber $F$ is complete with $\phi$-parabolic universal Riemannian covering, and such that the warping function $f$ is monotone. The only entire bounded solutions of the following modified $\phi$-maximal spacelike hypersurface equation*

$$
(\mathrm{E}) \quad \begin{cases} \mathrm{div}_\phi \left( \dfrac{Du}{f(u)\sqrt{f(u)^2 - |Du|_F^2}} \right) = -\dfrac{f'(u)}{\sqrt{f(u)^2 - |Du|_F^2}} \left( n + \dfrac{|Du|_F^2}{f(u)^2} \right) \\ |Du|_F \leq \alpha f(u), \end{cases}
$$

*where $0 < \alpha < 1$ is constant, are the constant functions $u = t_0$, with $f'(t_0) = 0$.*

*Proof.* We observe that the constraint on $|Du|_F$ assures the boundedness of the hyperbolic angle function $\Theta(u)$ of $S(u)$. Indeed, from (2.3), (2.5) and (5.2) we obtain that

$$(5.4) \qquad |\nabla h|^2 = \frac{|Du|_F^2}{f^2(u) - |Du|_F^2}.$$

Hence, using (2.5) and (5.4) we see that $|Du|_F \leq \alpha f(u)$ implies $\Theta(u) \geq \dfrac{-1}{\sqrt{1-\alpha^2}}$.

Moreover, since we are looking for bounded solutions for (E), taking $c = (1-\alpha^2)\inf_{S(u)} f^2(u)$ we can apply Proposition 1 to see that such a solution must be complete. Therefore, the result follows from Theorem 2. □

It is not difficult to see that we can reason as in the proof of the previous result and obtain non-parametric versions of all others theorems of Section 4. For instance, we have the following non-parametric version of Theorem 7.

**Theorem 9.** *Let $M_\phi$ be a spatially weighted steady state type spacetime, whose fiber $F$ is complete with $\phi$-parabolic universal Riemannian covering. The only entire bounded solutions of the equation*

$$(\mathbf{E}') \begin{cases} \mathrm{div}_\phi\left(\dfrac{Du}{e^u\sqrt{e^{2u}-|Du|_F^2}}\right) = -nH_\phi(u) - \dfrac{e^u}{\sqrt{e^{2u}-|Du|_F^2}}\left(n + \dfrac{|Du|_F^2}{e^{2u}}\right) \\ H_\phi(u) \geq 1 \\ |Du|_F \leq \sqrt{1 - \dfrac{1}{H_\phi^2(u)}}e^u \end{cases}$$

*are the constant functions $u = t_0$.*

*Proof.* Let us observe first that, under the assumptions of the theorem, the solutions of $(\mathbf{E}')$ determine complete entire graphs $S(u)$. In fact, from the last inequality of $(\mathbf{E}')$ we easily obtain that

$$e^{2u} - |Du|_F^2 \geq c = e^{2\inf_{S(u)} u} \sup_{S(u)} \dfrac{1}{H_\phi^2(u)} > 0,$$

and the completeness of $S(u)$ follows again from Proposition 1.

Therefore, using (5.4) and the last inequality of $(\mathbf{E}')$ we can easily verify that the hyperbolic angle function $\Theta(u)$ of $S(u)$ satisfies $-\Theta(u) \leq H_\phi(u)$ and, hence, the result follows from Theorem 7. □

## Acknowledgements

The first author is partially supported by MINECO-FEDER Grant No. MTM2012-34037. The second author is partially supported by CNPq, Brazil, grant 300769/2012-1. The third author is partially supported by CAPES, Brazil. The fourth author is partially supported by CAPES/CNPq, Brazil, grant Casadinho/Procad 552.464/2011-2.

[1] Departamento de Matemáticas, Campus Universitario de Rabanales, Universidad de Córdoba, 14071 Córdoba, Spain.

*E-mail address*: `alma.albujer@uco.es`

[2] Departamento de Matemática, Universidade Federal de Campina Grande, 58429-970 Campina Grande, Paraíba, Brazil.

*E-mail address*: `henrique@mat.ufcg.edu.br`
*E-mail address*: `arlandsonm@gmail.com`
*E-mail address*: `marco.velasquez@pq.cnpq.br`